\documentclass[10pt,twoside]{article}
\usepackage{amsmath,amssymb}
\usepackage[latin1]{inputenc}

\def\YEAR{\year}\newcount\VOL\VOL=\YEAR\advance\VOL by-1995
\def\firstpage{1}\def\lastpage{1000}
\def\received{}\def\revised{}
\def\communicated{}

\makeatletter
\def\magnification{\afterassignment\m@g\count@}
\def\m@g{\mag=\count@\hsize6.5truein\vsize8.9truein\dimen\footins8truein}
\makeatother

\font\eightrm=cmr8
\font\caps=cmcsc10                    
\font\Caps=cmcsc10 scaled \magstep1   

\makeatletter
\@specialpagefalse\headheight=8.5pt
\def\DocMath{}
\renewcommand{\@evenhead}{%
    \ifnum\thepage>\lastpage\rlap{\thepage}\hfill%
    \else\rlap{\thepage}\slshape\leftmark\hfill{\caps\SAuthor}\hfill\fi}%
\renewcommand{\@oddhead}{%
    \ifnum\thepage=\firstpage{\DocMath\hfill\llap{\thepage}}%
    \else{\slshape\rightmark}\hfill{\caps\STitle}\hfill\llap{\thepage}\fi}%
\makeatother

\def\TSkip{\bigskip}
\newbox\TheTitle{\obeylines\gdef\GetTitle #1
\ShortTitle  #2
\SubTitle    #3
\Author      #4
\ShortAuthor #5
\EndTitle
{\setbox\TheTitle=\vbox{\baselineskip=20pt\let\par=\cr\obeylines%
\halign{\centerline{\Caps##}\cr\noalign{\medskip}\cr#1\cr}}%
        \copy\TheTitle\TSkip\TSkip%
\def\next{#2}\ifx\next\empty\gdef\STitle{#1}\else\gdef\STitle{#2}\fi%
\def\next{#3}\ifx\next\empty%
    \else\setbox\TheTitle=\vbox{\baselineskip=20pt\let\par=\cr\obeylines%
    \halign{\centerline{\caps##} #3\cr}}\copy\TheTitle\TSkip\TSkip\fi%
\centerline{\caps #4}\TSkip\TSkip%
\def\next{#5}\ifx\next\empty\gdef\SAuthor{#4}\else\gdef\SAuthor{#5}\fi%
\ifx\received\empty\relax
    \else\centerline{\eightrm Received: \received}\fi%
\ifx\revised\empty\TSkip%
    \else\centerline{\eightrm Revised: \revised}\TSkip\fi%
\ifx\communicated\empty\relax
    \else\centerline{\eightrm Communicated by \communicated}\fi\TSkip\TSkip%
\catcode'015=5}}\def\Title{\obeylines\GetTitle}
\def\Abstract{\begingroup\narrower
    \parskip=\medskipamount\parindent=0pt{\caps Abstract. }}
\def\EndAbstract{\par\endgroup\TSkip}

\long\def\MSC#1\EndMSC{\def\arg{#1}\ifx\arg\empty\relax\else
     {\par\narrower\noindent%
     2000 Mathematics Subject Classification: #1\par}\fi}

\long\def\KEY#1\EndKEY{\def\arg{#1}\ifx\arg\empty\relax\else
        {\par\narrower\noindent Keywords and Phrases: #1\par}\fi\TSkip}

\newbox\TheAdd\def\Addresses{\vfill\copy\TheAdd\vfill
    \ifodd\number\lastpage\vfill\eject\phantom{.}\vfill\eject\fi}
{\obeylines\gdef\GetAddress #1
\Address #2 
\Address #3
\Address #4
\EndAddress
{\def\xs{4.3truecm}\parindent=0pt
\setbox0=\vtop{{\obeylines\hsize=\xs#1\par}}\def\next{#2}
\ifx\next\empty 
     \setbox\TheAdd=\hbox to\hsize{\hfill\copy0\hfill}
\else\setbox1=\vtop{{\obeylines\hsize=\xs#2\par}}\def\next{#3}
\ifx\next\empty 
     \setbox\TheAdd=\hbox to\hsize{\hfill\copy0\hfill\copy1\hfill}
\else\setbox2=\vtop{{\obeylines\hsize=\xs#3\par}}\def\next{#4}
\ifx\next\empty\ 
     \setbox\TheAdd=\vtop{\hbox to\hsize{\hfill\copy0\hfill\copy1\hfill}
                \vskip20pt\hbox to\hsize{\hfill\copy2\hfill}}
\else\setbox3=\vtop{{\obeylines\hsize=\xs#4\par}}
     \setbox\TheAdd=\vtop{\hbox to\hsize{\hfill\copy0\hfill\copy1\hfill}
                \vskip20pt\hbox to\hsize{\hfill\copy2\hfill\copy3\hfill}}
\fi\fi\fi\catcode'015=5}}\gdef\Address{\obeylines\GetAddress}

\begin{document}
\Title   Enumerative properties of generalized associahedra
\ShortTitle 
\SubTitle 
\Author  Frédéric Chapoton
\ShortAuthor  
\EndTitle
\Abstract 
Some enumerative aspects of the fans, called generalized associahedra,
introduced by S. Fomin and A. Zelevinsky in their theory of cluster
algebras are considered, in relation with a bicomplex and its two
spectral sequences. A precise enumerative relation with the lattices
of generalized noncrossing partitions is conjectured and some evidence
is given. \EndAbstract \MSC
\EndMSC
\KEY 
Generalized associahedra, noncrossing partition, f-vector
\EndKEY
\Address 
Institut Girard Desargues
Université Claude Bernard (Lyon 1)
21 Avenue Claude Bernard
F-69622 VILLEURBANNE Cedex
FRANCE
\Address
\Address
\Address
\EndAddress
\newcommand{\ZZ}{\mathbb{Z}}
\newcommand{\h}{\textrm{h}}
\newcommand{\coh}{\textrm{H}^*}
\newcommand{\none}{\textrm{NN}}
\newcommand{\gr}{\operatorname{gr}}
\newcommand{\zero}{\widehat{0}}
\newcommand{\one}{\widehat{1}}
\newcommand{\LW}{\textsf{L}_W}
\newcommand{\phip}{\Phi_{\geq -1}}
\newcommand{\rk}{\operatorname{rk}}

\newcommand{\dinde}[2]{\genfrac{}{}{0pt}{}{#1}{#2}}

\newtheorem{claim}{Claim}
\newtheorem{prop}{Proposition}
\newtheorem{theorem}{Theorem}
\newtheorem{conjecture}{Conjecture}
\newtheorem{lemma}{Lemma}
\newtheorem{coro}{Corollary}

\newenvironment{proof}{\begin{trivlist}\item{\bf{Proof.}}}
  {\hfill\rule{2mm}{2mm}\end{trivlist}}

\setcounter{section}{-1}

\section{Introduction}

In their work on cluster algebras \cite{cluster1,cluster2,ysystem}, S.
Fomin and A. Zelevinsky have introduced simplicial fans associated to
finite crystallographic root systems. These fans are associated with
convex polytopes called generalized associahedra \cite{polytopal} and
have been shown to be related to classical combinatorial objects such
as triangulations, noncrossing and nonnesting partitions and Catalan
numbers. The lattice of noncrossing partitions, which was defined
first for symmetric groups by G. Kreweras \cite{kreweras}, has been
recently generalized to all finite Coxeter groups
\cite{bessis,biane,bradywatt}. Surveys of its properties can be found
in \cite{mccammond} and \cite{simionNC}.

The aim of the present article is twofold. First, a refined
enumerative invariant, called the $F$-triangle, of the fan associated
to a root system is introduced and an inductive procedure is given for
its computation. The $F$-triangle is then related to a simple
combinatorial bicomplex and its spectral sequences. The second theme
is a conjecture which relates, through an explicit change of
variables, the $F$-triangle of a root system and a bivariate polynomial
defined in terms of the noncrossing partition lattice for the
corresponding Weyl group. Several evidences are given for this
conjecture.

The final section contains the computation of the $F$-triangle for the
root systems of type $A$ and $B$, using arguments based on
hypergeometric functions.

\smallskip

Thanks to C. Krattenthaler for his help with hypergeometric identities.

\section{The simplicial fans of clusters}

Let $\Phi$ be the root system associated to an irreducible Dynkin
diagram $X_n$ of finite type and rank $n$. Thus $X_n$ is among the
Killing-Cartan list $A_n,B_n,C_n,D_n$ or $E_6,E_7,E_8,F_4,G_2$. Let
$I$ be the underlying set of the Dynkin diagram and $\{\alpha_i\}_{i
  \in I}$ be the set of simple positive roots in $\Phi$.

Let us recall briefly the construction by S. Fomin and A. Zelevinsky
of the simplicial fan $\Delta(\Phi)$. Let $\phip$ be the union of the
set of negative simple roots $\{-\alpha_i\}_{i \in I}$ with the set
$\Phi_{>0}$ of positive roots. Elements of $\phip$ are called almost
positive roots. A symmetric binary relation on $\phip$ called
compatibility was defined in \cite{ysystem}. The following is
\cite[Theorem 1.10]{ysystem}.

\begin{theorem}
  The cones spanned by subsets of mutually compatible elements in
  $\phip$ define a complete simplicial fan $\Delta(\Phi)$.
\end{theorem}

From now on, cones of the fan $\Delta(\Phi)$ will be identified with
their spanning set of mutually compatible elements of $\phip$. The
cones of dimension $n$ of $\Delta(\Phi)$ are in bijection with maximal
mutually compatible subsets of $\phip$, which are called
\textit{clusters}. The cones of dimension $1$ of $\Delta(\Phi)$ are in
bijection with $\phip$ and will be called roots. A cone of
$\Delta(\Phi)$ is called \textit{positive} if it is spanned by
positive roots and \textit{non-positive} else.

One can define a fan $\Delta(\Phi)$ also for a non-irreducible root
system $\Phi$, as the product of the fans associated to its
irreducible components.

Let $P$ be the closed cone spanned by simple positive roots. This is
not a cone of $\Delta(\Phi)$ in general.

\begin{prop}
  \label{posi}
  The cone $P$ is exactly the union of all positive cones of $\Delta(\Phi)$.
\end{prop}
\begin{proof}
  Each positive cone is spanned by positive roots, hence is contained
  in $P$. Conversely, as the fan is complete, $P$ is contained in the
  union of all cones whose interior meet $P$. The interior of a
  non-positive cone does not meet $P$ as it consists of vectors with
  at least one negative coordinate in the basis of simple roots. Hence
  $P$ is contained in the union of all positive cones.
\end{proof}

As a special case of the description of the fan $\Delta(\Phi)$ in
\cite{ysystem}, the following Lemma holds.

\begin{lemma}
  \label{negaspan}
  The span of negative simple roots is a cone of $\Delta(\Phi)$
\end{lemma}

We recall \cite[Proposition 3.6]{ysystem} for later use.

\begin{prop}
  \label{supportnega}
  For every subset $J \subseteq I$, the correspondence $c \mapsto c
  \setminus \{-\alpha_i \}_{i \in J}$ is a bijection between cones of
  $\Delta(\Phi)$ whose negative part is $J$ and positive cones of
  $\Delta(\Phi(I\setminus J))$, where $\Phi(I\setminus J)$ is the
  restriction of the root system $\Phi$ to $I \setminus J$.
\end{prop}

\section{The $F$-triangle and the bicomplex of cones}

\label{spectral}

Let us define the $F$-triangle by its generating function
\begin{equation}
  F(\Phi)=F(x,y)=\sum_{k=0}^{n}\sum_{\ell=0}^{n}f_{k,\ell}x^k y^\ell,
\end{equation}
where $f_{k,\ell}$ is the cardinality of the set $C_{k,\ell}$ of cones
of $\Delta(\Phi)$ spanned by exactly $k$ positive roots and $\ell$
negative simple roots. The coefficient $f_{k,\ell}$ vanishes if
$k+\ell > n$, hence the name triangle.

\begin{prop}
  \label{recuF}
  The $F$-triangle has the following properties.
  \begin{enumerate}
  \item If $\Phi$ and $\Phi'$ are two root systems, one has $F(\Phi
    \times \Phi')=F(\Phi)\times F(\Phi')$.
  \item If $\Phi$ is an irreducible root system on $I$, then one has
    \begin{equation}
      \partial_y F (\Phi(I))=\sum_{i \in I} F(\Phi(I\setminus \{i\})),
    \end{equation}
    where $\Phi(I\setminus \{i\})$ is the restriction of the root
    system $\Phi$ to $I\setminus \{i\}$.
  \end{enumerate}
\end{prop}
\begin{proof}
  The first statement is obvious. The proof of the second statement is
  by double counting of the sets of pairs
  \begin{equation}
    \{ (i,c) \mid -\alpha_i \in c \text{ and } c \in C_{k,\ell}\},
  \end{equation}
  for all $k$ and $\ell$. Let us fix $k$ and $\ell$.
  
  On one side, the cardinalityof this set is just $\ell f_{k,\ell}$ by
  definition of $C_{k,\ell}$. On the other hand, by \cite[Proposition
  3.5 (3)]{ysystem}, the cardinalityis given by the sum over $i\in I$ of
  $f^i_{k,\ell-1}$ where $f^i$ is the $F$-triangle for the root system
  induced on $I \setminus \{i\}$.

  This gives the equality
  \begin{equation}
    \ell f_{k,\ell}= \sum_{i \in I} f^i_{k, \ell-1},
  \end{equation}
  which proves the second assertion of the proposition.
\end{proof}

\smallskip

The usual $f$-vector is given by the generating series
\begin{equation}
  \label{diagonal}
  f(x)=\sum_{k=0}^{n}f_{k}x^k=F(x,x),
\end{equation}
where $f_k$ is the number of cones of dimension $k$.

The following is \cite[Proposition 3.7]{ysystem}.

\begin{prop}
  \label{recuf}
  The $f$-vector has the following properties.
  \begin{enumerate}
  \item If $\Phi$ and $\Phi'$ are two root systems, one has $f(\Phi
    \times \Phi')=f(\Phi)\times f(\Phi')$.
  \item If $\Phi$ is an irreducible root system on $I$, then one has
    \begin{equation}
      \partial_x f (\Phi(I))=\frac{\h+2}{2}
      \sum_{i \in I} f(\Phi(I\setminus \{i\})),
    \end{equation}
    where $\h$ is the Coxeter number of $\Phi$.
  \end{enumerate}
\end{prop}

\smallskip

Together Proposition \ref{recuF} and Proposition \ref{recuf} are
sufficient to compute simultaneously the $F$-triangle and the
$f$-vector by induction on the cardinality of $I$, using
(\ref{diagonal}).

For example, the $A_3$ $f$-vector is $(1,9,21,14)$ and the $A_3$
$F$-triangle is presented below, with $k=0$ to $3$ from top to bottom
and $\ell=0$ to $3$ from left to right.
\begin{equation}
  \begin{bmatrix}
    1&3&3&1\\
    6&8&3\\
    10&5\\
    5
  \end{bmatrix}
\end{equation}

\smallskip

The $F$-triangle has a nice symmetry property, which is a refined
version of the classical Dehn-Sommerville equations for complete
simplicial fans.

\begin{prop}
  \label{principale}
  One has 
  \begin{equation}
    F(x,y)=(-1)^n F(-1-x,-1-y).
  \end{equation}
\end{prop}
\begin{proof}
  The proof is just an adaptation of the original proof of the
  Dehn-Sommerville equations (see \cite[p. 212-213]{babil}), taking
  care of the two different kind of half-edges of the fan. The
  proposition is equivalent to the set of equations
  \begin{equation}
    f_{i,j}=\sum_{k,\ell} (-1)^{n+k+\ell}\binom{k}{i}\binom{\ell}{j}f_{k,\ell},
  \end{equation}
  for all $i,j$. Let us now fix $i$ and $j$ and compute
  $f_{i,j}$. First, it is given by
  \begin{equation}
    f_{i,j}=\sum_{c\in C_{i,j}} 1.
  \end{equation}
  Then using Lemma \ref{complet}, this is rewritten as
  \begin{equation}
    \sum_{c\in C_{i,j}} \sum_{c \subseteq d} (-1)^{n+\dim(d)}.
  \end{equation}
  Exchanging summations, one obtains
  \begin{equation}
    \sum_{k,\ell} (-1)^{n+k+\ell} \sum_{d \in C_{k,\ell}}
    \sum_{\dinde{c\subseteq d}{c \in C_{i,j}}} 1.
  \end{equation}
  As the fan is simplicial, this is
  \begin{equation}
    \sum_{k,\ell} (-1)^{n+k+\ell} \binom{k}{i}\binom{\ell}{j}f_{k,\ell}.
  \end{equation}
  The proposition is proved.
\end{proof}

The following lemma is classical. It follows from the fact that the
link of a simplex in a homology sphere is again a homology sphere, see
\cite[p. 214]{babil}.

\begin{lemma}
  \label{complet}
  Let $c$ be a cone in a complete simplicial fan of dimension $n$.
  Then
  \begin{equation}
    \sum_{c \subseteq d}(-1)^{n+\dim(d)}=1.
  \end{equation}
\end{lemma}

Let us now introduce two specializations of the $F$-triangle.

The positive $f^+$-vector is given by the generating series
\begin{equation}
  f^+(x)=\sum_{k=0}^{n}f^+_{k}x^k=F(x,0),
\end{equation}
where $f^+_k=f_{k,0}$ is the number of positive cones of dimension
$k$.

The natural $f^\natural$-vector is given by the generating series
\begin{equation}
  f^\natural(x)=\sum_{k=0}^{n}f^\natural_{k}x^k=F(x,-1).
\end{equation}
Its interpretation will be given later in Proposition \ref{natural}.

\smallskip

The symmetry obtained in Proposition \ref{principale} has the following
consequences.

\begin{coro}
  The $f^+$-vector and $f^\natural$-vector determine each other. One has
 \begin{equation}
    F(x,0)=(-1)^n F(-1-x,-1).
  \end{equation}
\end{coro}

\begin{coro}
  \label{coro2}
  One has
  \begin{equation}
    F(0,x)=(-1)^n F(-1,-1-x)=(x+1)^n.
  \end{equation}
  Equivalently, one has $F(-1,y)=y^n$.
\end{coro}
\begin{proof}
  By Lemma \ref{negaspan}, each subset of the set of negative simple roots is
  a cone of $\Delta(\Phi)$.
\end{proof}


\smallskip

Let us define a bicomplex on cones and study its two spectral
sequences. This bicomplex is essentially the complex associated to the
simplicial set defined by the fan, where the differential is split
according to the two kinds of elements of $\phip$.

In the unital exterior algebra generated over $\ZZ$ by the set
$\phip$, consider the linear span $D$ of the monomials associated to
the cones of the fan $\Delta(\Phi)$. The simplicial differential of a
monomial in $D$ is defined by
\begin{equation}
  d(\alpha_1 \wedge \alpha_2 \wedge \dots \wedge \alpha_k)=
  \sum_{\ell=1}^k (-1)^{\ell-1} 
  \alpha_1 \wedge \alpha_2 \wedge \dots \widehat{\alpha_\ell}
 \dots \wedge \alpha_k,
\end{equation}
where $\widehat{\alpha_\ell}$ means that $\alpha_\ell$ has been
removed.

This gives a complex $(D,d)$ computing the reduced homology of a
sphere.

The differential $d$ is the sum of two maps $d_+$ and $d_-$ which
correspond respectively to the removal of a positive root or a
negative simple root in a monomial. It is clear that each of $d_+$ and
$d_-$ is a differential. This defines a bicomplex structure on $D$
when bigraded by the number of positive roots and negative simple
roots.

\begin{prop}
  \label{primo}
  The spectral sequence of $D$ starting with $d_+$ degenerates at
  first step.
\end{prop}
\begin{proof}
  The complex $(D,d_+)$ decomposes as a direct sum of complexes
  $D_{S_-}$ according to the fixed negative part $S_-$. If this
  negative part $S_-$ is not the full set $\{-\alpha_i\}_{i \in I}$,
  the subcomplex $D_{S_-}$ has no homology. To prove this, it is
  enough to consider the complex of positive cones for the
  differential $d_+$ for all root systems, as the subcomplex $D_{S_-}$
  with fixed negative part $S_-$ is isomorphic to the complex of
  positive cones for a smaller root system by Proposition
  \ref{supportnega}. The complex of positive cones for the
  differential $d_+$ has no homology because it computes the reduced
  homology of the contractible simplicial complex $P$, by Proposition
  \ref{posi}.
\end{proof}

\begin{prop}
  \label{secundo}
  The spectral sequence of $D$ starting with $d_-$ degenerates at
  second step.
\end{prop}
\begin{proof}
  The complex $(D,d_-)$ decomposes as a direct sum of complexes
  $D_{S_+}$ according to the fixed positive part $S_+$ of pairwise
  compatible positive roots. Let $N(S_+)$ be the set of negative
  simple roots compatible with all roots of $S_+$. Then the complex
  $D_{S_+}$ is isomorphic to the tensor product over the set $N(S_+)$
  of contractible complexes of the following shape
  \begin{equation}
    0 \longrightarrow \ZZ \stackrel{\simeq}{\longrightarrow}
    \ZZ \longrightarrow 0.
  \end{equation}
  Hence either $N(S_+)$ is empty and $d_-$ vanishes on $D_{S_+}$, or
  the subcomplex $D_{S_+}$ has no homology. Therefore the homology of
  $(D,d_-)$ is concentrated in positive cones and is given exactly by
  positive cones which can not be extended with negative simple roots.
  This implies the collapsing of the spectral sequence at second step.
\end{proof}

A positive cone of which is not a subcone of a non-positive cone is
called a \textit{natural cone}.

\begin{prop}
  \label{natural}
  The number of natural cones of dimension $k$ is $f^\natural_k$.
\end{prop}
\begin{proof}
  The natural $f^\natural$-vector is by definition the Euler
  characteristic of the complex $(D,d_-)$. By the proof of Proposition
  \ref{secundo}, the homology of $d_-$ is concentrated in degree $0$
  and has dimension given by the numbers of natural cones.
\end{proof}

The first step of the spectral sequence starting with $d_-$ is
therefore a complex on natural cones. Its homology is concentrated in
degree $n$.

\section{The Lattice of noncrossing partitions}

The lattice $\LW$ of noncrossing partitions associated to a finite
Coxeter group $W$ has been introduced independently in
\cite{bessis,biane,bradywatt}. The surveys \cite{mccammond,simionNC}
give a good feeling of its importance in different parts of
mathematics. Let us recall shortly its definition.

Let $n$ be the rank of $W$. Let $S=\{s_1,\dots,s_n\}$ be the set of
simple reflections in $W$. Let $T$ be the set of all reflections in
$W$. As $W$ is also generated by $T$, one can define a length function
$\ell_T$ on $W$ with respect to the generators in $T$. Using this
length function, a partial order is defined on $W$ as follows. Let $v$
and $w$ be two elements of $W$. Then $v \leq w$ if and only if there
exists $t_1, \dots ,t_k$ in $T$ such that $w=t_1 \dots t_k v$ and
$\ell_T(w)=\ell_T(v)+k$. Maximal elements for this partial order are
exactly Coxeter elements of $W$, \textit{i.e.} products in some order
of the set $S$ of simple reflections. The group $W$ acting by
conjugation gives automorphisms of this partial order, as the set $T$
of all reflections is stable by conjugation. As all Coxeter elements
are conjugated, one can define, up to isomorphism, the noncrossing
partition lattice $\LW$ to be the interval between the unit and a
Coxeter element for this partial order.

From this description, one can see that the noncrossing partition
poset associated to the product of two Coxeter groups is isomorphic to
the product of the noncrossing partition posets associated to these
groups.

We shall now recall some properties of $\LW$. The poset $\LW$ is a
finite lattice which is graded of rank $n$ and self-dual. Let $\zero$
and $\one$ be the minimum and maximum elements of $\LW$. Let $\mu$ be
the Möbius function of $\LW$. Then the lattice $\LW$ has the following
invariants.

\begin{prop}
  \label{invariants}
  The Zeta polynomial of $\LW$ is
  \begin{equation}
    Z_W(X) =\prod_{i =1}^{n}\frac{\h X-e_i+1}{e_i+1},
  \end{equation}
  the cardinalityof $\LW$ is 
  \begin{equation}
    \# \LW = \prod_{i =1}^{n}\frac{\h+e_i+1}{e_i+1},
  \end{equation}
  and the Möbius number of $\LW$ is
  \begin{equation}
    \label{mobiusnr}
    \mu(\zero,\one) = (-1)^n \prod_{i =1}^{n}\frac{\h+e_i-1}{e_i+1},
  \end{equation}
  where $\h$ is the Coxeter number and $e_1,\dots,e_n$ are the
  exponents of $W$.
\end{prop}

This proposition follows from the known formulas for the Zeta
polynomial of $\LW$ in the classical cases \cite{athareiner,reiner}
and from the computation of the Zeta polynomials in the exceptional
cases by V. Reiner \cite{reinerprivate}. The statements about cardinal
and Möbius numbers follows from the knowledge of the Zeta polynomial.
The cardinality part of the Proposition can also be checked from the
data in \cite{picantin}. To find a uniform proof of Proposition
\ref{invariants} is an interesting open problem.

Let $\rk$ be the rank function of $\LW$. Consider the following
generating function for Möbius numbers of intervals in $\LW$ according
to the ranks:
\begin{equation}
  M(x,y)=\sum_{a \leq b} \mu(a,b) x^{\rk(b)} y^{\rk(a)}.
\end{equation}
This generating function is called the $M$-triangle for $W$.

Let us assume from now on that $W$ is the Weyl group of a
crystallographic root system $\Phi$. Here is the main Conjecture.

\begin{conjecture}
  \label{cluster-treillis}
  The $F$-triangle for $\Phi$ and the $M$-triangle for $W$ are related
  by the following invertible transformation:
  \begin{equation}
    \label{mystere}
    (1-y)^n F(\frac{x+y}{1-y},\frac{y}{1-y})=M(-x,-y/x).
  \end{equation}
\end{conjecture}

Let us note that the left-hand side of (\ref{mystere}) can be
rewritten using Proposition \ref{principale} as
\begin{equation}
  \label{alternative}
  (y-1)^n F(\frac{x+1}{y-1},\frac{1}{y-1}).
\end{equation}

It is not hard to check by hand that Conjecture \ref{cluster-treillis}
holds for root systems of small ranks. It is probably possible to
prove it for classical types using the combinatorial descriptions of
the noncrossing lattices \cite{athareiner,reiner} and of the
generalized associahedra \cite{ysystem}. Using a computer, one could
check it for most of the exceptional types. Rather than doing that, we
prefer to give now several conceptual evidences for Conjecture
\ref{cluster-treillis}.

\subsection{First Evidence}

Let us consider the value of (\ref{mystere}) at $x=0$. The left-hand
side becomes $(1-y)^n F(\frac{y}{1-y},\frac{y}{1-y})$, which is
nothing but the $h$-vector of the simplicial fan $\Delta(\Phi)$. The
right-hand side is
\begin{equation}
  \sum_{a} y^{\rk(a)},
\end{equation}
which is known to coincide with the $h$-vector for all root systems,
see \cite{athanams,athareiner,bessis}. These $h$-vectors are sometimes
called the generalized Narayana numbers. Note that there is no uniform
proof of this fact known so far.

\subsection{Second evidence}

Let us compute the coefficient of $x^n$ in the constant term of
(\ref{mystere}) with respect to $y$. On the left-hand side, this is
the number $f_{n,0}$ of positive clusters, which is known by
\cite[Proposition 3.9]{ysystem} to be given by
\begin{equation}
  \prod_{i =1}^{n}\frac{\h+e_i-1}{e_i+1}.
\end{equation}
On the right-hand-side, one gets $(-1)^n$ times the Möbius number
$\mu(\zero,\one)$, which is given by (\ref{mobiusnr}).

\subsection{Third evidence}

It is known that the lattice $\LW$ is self-dual \cite[\S 2.3]{bessis}.
This fact implies the following symmetry of the $M$-triangle:
\begin{equation}
  M(x,y)=(xy)^n M(1/y,1/x).
\end{equation}

\begin{claim}
  Assuming Conjecture \ref{cluster-treillis}, this symmetry is
  equivalent to Proposition \ref{principale}.
\end{claim}

Indeed, the symmetry of the $M$-triangle is equivalent to
\begin{equation}
  M(-x,-y/x)=y^n M(-x/y,-1/x).
\end{equation}

Now applying Conjecture \ref{cluster-treillis} to the right-hand side,
one gets (\ref{alternative}), which is the same as (\ref{mystere})
thanks to Proposition \ref{principale}. This proves the claim.

\subsection{Fourth evidence}

Consider the value of (\ref{mystere}) at $x=-1$. The left-hand side is
computed using (\ref{alternative}) to be
\begin{equation}
  (y-1)^n F(0,\frac{1}{y-1}).
\end{equation}
But this is just $y^n$ by Corollary \ref{coro2}. That the right-hand
side is also $y^n$ is a standard property of the $M$-triangle for
graded posets of rank $n$ with $\zero$ and $\one$.

\subsection{Fifth evidence}

By Proposition \ref{recuF}, the $F$-triangle has a multiplicative
behavior with respect to product of root systems. It is classical
that the $M$-triangle is also multiplicative for the product of
lattices. These multiplicative behaviors are compatible with Formula
(\ref{mystere}).

\section{Computation for type $A$}

In this section, the $F$-triangle is computed for root systems of type
$A$. This can serve as a first step towards the proof of Conjecture
\ref{cluster-treillis} in type $A$.

Let us first recall the known expression for the $f$-vector in type
$A$ \cite{lee,simionB}.

\begin{prop}
  The $f$-vector for $A_n$ is given by
  \begin{equation}
    \sum_{k=0}^{n} \frac{1}{k+1}\binom{n}{k}\binom{n+k+2}{k}x^k.
  \end{equation}
\end{prop}

Let us define a generating function for $f$-vectors of type $A$. First
the $f$-vector for $A_n$ is made homogeneous of degree $n$ using a new
variable $z$, then all homogenized $f$-vectors are added. Let
\begin{equation}
  f=\sum_{m=0}^{\infty} \sum_{k=0}^{\infty} 
  \frac{1}{k+1}\binom{m+k}{k}\binom{2k+m+2}{k}x^k z^m.
\end{equation}

Our aim is now to prove the following Proposition.

\begin{prop}
  \label{FvectorA}
  The $F$-triangle for $A_n$ is given by
  \begin{equation}
    \label{FtypeA}
    \sum_{k=0}^{n}\sum_{\ell=0}^{n} 
    \frac{\ell+1}{k+\ell+1}\binom{n}{k+\ell}\binom{n+k}{n}x^k y^\ell.
  \end{equation}
\end{prop}

Let us define similarly a generating function for the functions
(\ref{FtypeA}). Let
\begin{equation}
  F=\sum_{\ell=0}^{\infty} \sum_{m=0}^{\infty} \sum_{k=0}^{\infty} 
  \frac{\ell+1}{k+\ell+1}\binom{k+\ell+m}{k+\ell}\binom{2k+\ell+m}{k+\ell+m}
  x^k y^\ell z^m.
\end{equation}

Recall that Proposition \ref{recuF} and Formula (\ref{diagonal}) give
a recursion for computing the $F$-triangle assuming that the
$f$-vector is known. In type $A$, the induction given by Proposition
\ref{recuF} is easily seen to be equivalent to the equation
$\partial_y F=F^2$. Hence, to prove Proposition \ref{FvectorA}, it is
enough to prove that $\partial_y F=F^2$ and that the substitution of
$y$ by $x$ in $F$ is $f$. This is done in the next two Lemmas.

\begin{lemma}
  The partial derivative $\partial_y F$ equals $F^2$.
\end{lemma}
\begin{proof}
  Let us fix $k,\ell,m$ and consider the coefficient of $x^k y^\ell
  z^m$ in $F^2$. It is given by
  \begin{equation}
    \sum_{k_1,\ell_1,m_1}
    \frac{(\ell_1+1)(2k_1+\ell_1+m_1)!}{(k_1+\ell_1+1)!m_1!k_1!}
    \frac{(\ell-\ell_1+1)(2k-2k_1+\ell-\ell_1+m-m_1)!}
    {(k-k_1+\ell-\ell_1+1)!(m-m_1)!(k-k_1)!}.
  \end{equation}
  One can first do the summation with respect to $m_1$ using the
  Chu-Vandermonde identity. The result is
  \begin{multline}
    \binom{2k+\ell+m+1}{m}\sum_{k_1,\ell_1}
    \frac{\ell_1+1}{2k_1+\ell_1+1}\binom{2k_1+\ell_1+1}{k_1}\\
    \frac{\ell-\ell_1+1}{2k-2k_1+\ell-\ell_1+1}
    \binom{2k-2k_1+\ell-\ell_1+1}{k-k_1}.
  \end{multline}
  Then using a result of Carlitz \cite[Theorem 6, Formula
  (5.14)]{carlitz} with parameters
  $a=2,c=1,\alpha=\alpha'=1,b=0,d=0,\beta=\beta'=-1$, the remaining
  double sum can be computed. The result is
  \begin{equation}
    \binom{2k+\ell+m+1}{m}(-1)^\ell \frac{-2\ell-4}{(2k+\ell+2)(-2)}
    \binom{2k+\ell+2}{k}\binom{-2}{\ell}.
  \end{equation}
  This can be rewritten as
  \begin{equation}
    \frac{(2k+\ell+m+1)!(\ell+2)(\ell+1)}{m! k! (k+\ell+2)!},
  \end{equation}
  which is exactly the coefficient of $x^k y^\ell z^m$ in $\partial_y
  F$. The Lemma is proved.
\end{proof}

\begin{lemma}
  The substitution of $y$ by $x$ in $F$ equals $f$.
\end{lemma}
\begin{proof}
  Let us fix $K,m$ and compute the coefficient of $x^K z^m$ in this
  substituted $F$. This is given by
  \begin{equation}
    \frac{1}{(K+1)!m!}\sum_{\ell}\frac{(\ell+1)(2K-\ell+m)!}
    {(K-\ell)!}.
  \end{equation}
  This is rewritten as
  \begin{equation}
    \frac{(K+m)!}{(K+1)!m!}\sum_{\ell}\binom{\ell+1}{\ell}
    \binom{2K+m-\ell}{K-\ell}.
  \end{equation}
  Using the Chu-Vandermonde identity, this equals
  \begin{equation}
    \frac{(K+m)!}{(K+1)!m!}\binom{2K+m+2}{K},
  \end{equation}
  which in turn is equal to
  \begin{equation}
    \frac{(K+m)!(2K+m+2)!}{m!K!(K+1)!(K+m+2)!}.
  \end{equation}
  This is exactly the coefficient of $x^K y^\ell$ in $f$. The Lemma is
  proved.
\end{proof}

\section{Computation for type $B$}

In this section, the $F$-triangle is computed for root systems of type
$B$. Let us first recall the known expression for the $f$-vector in
type $B$ \cite{simionB}.

\begin{prop}
  The $f$-vector for $B_n$ is given by
  \begin{equation}
    \sum_{k=0}^{n} \binom{n}{k}\binom{n+k}{k}x^k.
  \end{equation}
\end{prop}

Let us define a generating function for $f$-vectors of type $B$ as we
did before for type $A$. By convention, $B_0$ is $A_0$ and $B_1$ is
$A_1$. Let
\begin{equation}
  g=\sum_{m=0}^{\infty} \sum_{k=0}^{\infty} 
  \binom{m+k}{k}\binom{2k+m}{k}x^k z^m.
\end{equation}

Our aim is now to prove the following Proposition.
\begin{prop}
  \label{FvectorB}
  The $F$-triangle for $A_n$ is given by
  \begin{equation}
    \label{FtypeB}
    \sum_{k=0}^{n}\sum_{\ell=0}^{n} 
    \binom{n}{k+\ell}\binom{n+k-1}{n-1}x^k y^\ell.
  \end{equation}
\end{prop}

Let us define similarly a generating function for the functions
(\ref{FtypeB}). Let
\begin{equation}
  G=\sum_{\ell=0}^{\infty} \sum_{m=0}^{\infty} \sum_{k=0}^{\infty} 
  \binom{k+\ell+m}{k+\ell}\binom{2k+\ell+m-1}{k+\ell+m-1}
  x^k y^\ell z^m.
\end{equation}

Recall again that Proposition \ref{recuF} and Formula (\ref{diagonal})
give a recursion for computing the $F$-triangle assuming that the
$f$-vector is known. In type $B$, the induction given by Proposition
\ref{recuF} is easily seen to be equivalent to the equation
$\partial_y G=F G$. Hence, to prove Proposition \ref{FvectorB}, it is
enough to prove that $\partial_y G=F G$ and that the substitution of
$y$ by $x$ in $G$ is $g$. This is done in the next two Lemmas.

\begin{lemma}
  The partial derivative $\partial_y G$ equals $F G$.
\end{lemma}
\begin{proof}
  The proof is very similar to the type $A$ case and will be more
  sketchy. Let us fix $k,\ell,m$ and consider the coefficient of $x^k
  y^\ell z^m$ in $F G$. It is given by
  \begin{multline}
    \label{fullB}
    \sum_{k_1,\ell_1,m_1}
    \frac{(\ell_1+1)(2k_1+\ell_1+m_1)!}{(k_1+\ell_1+1)!m_1!k_1!}\\
    \frac{(k-k_1+\ell-\ell_1+m-m_1)(2k-2k_1+\ell-\ell_1+m-m_1-1)!}
    {(k-k_1+\ell-\ell_1)!(m-m_1)!(k-k_1)!}.
  \end{multline}
  Let us consider separately the summation of factors depending on
  $m_1$ with respect to $m_1$:
  \begin{multline}
    \sum_{m_1}(k-k_1+\ell-\ell_1+m-m_1)(2k_1+\ell_1+m_1)!\\
    \frac{(2k-2k_1+\ell-\ell_1+m-m_1-1)!}
    {m_1!(m-m_1)!}.
  \end{multline}
  This sum can be split in two terms:
  \begin{multline}
    \sum_{m_1}\frac
    {(k-k_1+\ell-\ell_1)(2k_1+\ell_1+m_1)!(2k-2k_1+\ell-\ell_1+m-m_1-1)!}
    {m_1!(m-m_1)!}\\
    +    \sum_{m_1}\frac
    {(2k_1+\ell_1+m_1)!(2k-2k_1+\ell-\ell_1+m-m_1-1)!}
    {m_1!(m-m_1-1)!}.
  \end{multline}
  Rewriting the summations with binomial coefficients gives
  \begin{multline}
    (k-k_1+\ell-\ell_1)(2k_1+\ell_1)!(2k-2k_1+\ell-\ell_1-1)!
    \\ \sum_{m_1}
    \binom{2k_1+\ell_1+m_1}{m_1}\binom{2k-2k_1+\ell-\ell_1+m-m_1-1}
    {m-m_1}\\
    +    (2k_1+\ell_1)!(2k-2k_1+\ell-\ell_1)!\\ \sum_{m_1}
    \binom{2k_1+\ell_1+m_1}{m_1}\binom{2k-2k_1+\ell-\ell_1+m-m_1-1}{m-m_1-1}.
  \end{multline}
  Using the Chu-Vandermonde identity for each of these two terms gives
  \begin{multline}
    \label{sommem1}
    (k-k_1+\ell-\ell_1)(2k_1+\ell_1)!(2k-2k_1+\ell-\ell_1-1)!
    \binom{2k+\ell+m}{m}\\
    + (2k_1+\ell_1)!(2k-2k_1+\ell-\ell_1)!
    \binom{2k+\ell+m}{m-1}.
  \end{multline}
  Now it is time to plug this result into the full summation
  (\ref{fullB}). Let us do this separately for the the two terms of
  (\ref{sommem1}). The first term of (\ref{sommem1}) plugged into
  (\ref{fullB}) gives
  \begin{equation}
    \binom{2k+\ell+m}{m} \sum_{k_1,\ell_1}
    \binom{2(k-k_1)+(\ell-\ell_1)-1}{k-k_1}
    \frac{\ell_1+1}{2k_1+\ell_1+1}\binom{2k_1+\ell_1+1}{k_1}.
  \end{equation}
  Using another formula of Carltiz \cite[(5.15)]{carlitz} with
  parameters $a=2,c=1,b=0,d=0,\alpha=0,\beta=0,\alpha'=1,\beta'=-1$,
  this becomes
  \begin{equation}
    \label{term1}
    \binom{2k+\ell+m}{m} \binom{2k+\ell}{k}(\ell+1).
  \end{equation}
  The second term of (\ref{sommem1}) plugged into (\ref{fullB}) gives
  \begin{equation}
    \binom{2k+\ell+m}{m-1}
    \sum_{k_1,\ell_1}
    \binom{2(k-k_1)+(\ell-\ell_1)}{k-k_1}
    \frac{\ell_1+1}{2k_1+\ell_1+1}\binom{2k_1+\ell_1+1}{k_1}.
  \end{equation}
  Using again \cite[(5.15)]{carlitz} with parameters
  $a=2,c=1,b=0,d=0,\alpha=1,\beta=0,\alpha'=1,\beta'=-1$, this becomes
  \begin{equation}
    \label{term2}
    \binom{2k+\ell+m}{m-1}
    \binom{2k+\ell+1}{k}(\ell+1).
  \end{equation}
  Summing (\ref{term1}) and (\ref{term2}) gives
  \begin{equation}
    \frac{(\ell+1)(2k+\ell+m)!(k+\ell+m+1)}{k!m!(k+\ell+1)!}.
  \end{equation}
  This is exactly the coefficient of $x^k y^\ell z^m$ in $\partial_y
  G$. The Lemma is proved.
\end{proof}

\begin{lemma}
  The substitution of $y$ by $x$ in $G$ equals $g$.
\end{lemma}
\begin{proof}
  Let us fix $K,m$ and compute the coefficient of $x^K z^m$ in this
  substituted $G$. This is given by
  \begin{equation}
    \binom{K+m}{K}\sum_{k=0}^{K} \binom{K+k+m-1}{K+m-1}.
  \end{equation}
  This is rewritten using the standard column summation property of
  binomial coefficient as
  \begin{equation}
    \binom{K+m}{K} \binom{2K+m}{K+m}.
  \end{equation}
  This is exactly the coefficient of $x^K y^\ell$ in $g$. The Lemma is
  proved.
\end{proof}

\bibliographystyle{plain}
\bibliography{preuves}

\begin{thebibliography}{10}

\bibitem{athanams}
C.~Athanasiadis.
\newblock On a refinement of the generalized {C}atalan numbers for {W}eyl
  groups.
\newblock {\em T.A.M.S.}, 2004.

\bibitem{athareiner}
C.~Athanasiadis and V.~Reiner.
\newblock Noncrossing partitions for the group {$D_n$}.
\newblock 24 pages, 2003.

\bibitem{babil}
M.~M. Bayer and L.~J. Billera.
\newblock Counting faces and chains in polytopes and posets.
\newblock In {\em Combinatorics and algebra (Boulder, Colo., 1983)}, volume~34
  of {\em Contemp. Math.}, pages 207--252. Amer. Math. Soc., Providence, RI,
  1984.

\bibitem{bessis}
D.~Bessis.
\newblock The dual braid monoid.
\newblock {\em Ann. Sci. E.N.S.}, 2003.

\bibitem{biane}
P.~Biane.
\newblock Some properties of crossings and partitions.
\newblock {\em Discrete Math.}, 175(1-3):41--53, 1997.

\bibitem{bradywatt}
T.~Brady and C.~Watt.
\newblock {$K(\pi,1)$}'s for {A}rtin groups of finite type.
\newblock In {\em Proceedings of the Conference on Geometric and Combinatorial
  Group Theory, Part I (Haifa, 2000)}, volume~94, pages 225--250, 2002.

\bibitem{carlitz}
L.~Carlitz.
\newblock Some expansions and convolution formulas related to {M}ac{M}ahon's
  master theorem.
\newblock {\em SIAM J. Math. Anal.}, 8(2):320--336, 1977.

\bibitem{polytopal}
F.~Chapoton, S.~Fomin, and A.~Zelevinsky.
\newblock Polytopal realizations of generalized associahedra.
\newblock {\em Canad. Math. Bull.}, 45(4):537--566, 2002.
\newblock Dedicated to Robert V.\ Moody.

\bibitem{cluster1}
S.~Fomin and A.~Zelevinsky.
\newblock Cluster algebras. {I}. {F}oundations.
\newblock {\em J. Amer. Math. Soc.}, 15(2):497--529 (electronic), 2002.

\bibitem{cluster2}
S.~Fomin and A.~Zelevinsky.
\newblock Cluster algebras. {II}. {F}inite type classification.
\newblock {\em Inventiones Mathematicae}, 154:63--121, 2003.

\bibitem{ysystem}
S.~Fomin and A.~Zelevinsky.
\newblock Y-systems and generalized associahedra.
\newblock {\em Annals of Math.}, 158(3):977--1018, 2003.

\bibitem{kreweras}
G.~Kreweras.
\newblock Sur les partitions non crois\'ees d'un cycle.
\newblock {\em Discrete Math.}, 1(4):333--350, 1972.

\bibitem{lee}
C.~W. Lee.
\newblock The associahedron and triangulations of the {$n$}-gon.
\newblock {\em European J. Combin.}, 10(6):551--560, 1989.

\bibitem{mccammond}
J.~McCammond.
\newblock Noncrossing partitions in surprising locations.
\newblock 14 pages, 2003.

\bibitem{picantin}
M.~Picantin.
\newblock Explicit presentations for the dual braid monoids.
\newblock {\em C. R. Math. Acad. Sci. Paris}, 334(10):843--848, 2002.

\bibitem{reinerprivate}
V.~Reiner.
\newblock private communication.

\bibitem{reiner}
V.~Reiner.
\newblock Non-crossing partitions for classical reflection groups.
\newblock {\em Discrete Math.}, 177(1-3):195--222, 1997.

\bibitem{simionNC}
R.~Simion.
\newblock Noncrossing partitions.
\newblock {\em Discrete Math.}, 217(1-3):367--409, 2000.
\newblock Formal power series and algebraic combinatorics (Vienna, 1997).

\bibitem{simionB}
R.~Simion.
\newblock A type-{B} associahedron.
\newblock {\em Adv. in Appl. Math.}, 30(1-2):2--25, 2003.
\newblock Formal power series and algebraic combinatorics (Scottsdale, AZ,
  2001).

\end{thebibliography}
\end{document}